\documentclass[11pt]{article}
\title{Interleaved adjoints of directed graphs}

\author{Jan Foniok
\\ ETH Zurich, Institute for Operations Research
\\ R\"amistrasse 101, 8092 Zurich, Switzerland
\\ {\small Email: foniok@math.ethz.ch}
\and
Jaroslav Ne\v set\v ril\thanks{Partially
supported by the Project LN00A056 of the Czech Ministery 
of Education and by CRM Barcelona, Spain.}
\\ Department of Applied Mathematics
\\and
\\ Institute for Theoretical
Computer Science (ITI)
\\ Charles University
\\Malostransk\' e n\' am.25, 11800
Praha 1
\\ Czech Republic
\\{\small Email: nesetril@kam.ms.mff.cuni.cz}
\and
Claude Tardif\thanks{Supported by grants from 
NSERC and ARP}
\\
Department of Mathematics and Computer Science\\
       Royal Military College of Canada\\
       PO Box 17000, Station ``Forces''\\
       Kingston, Ontario K7K 7B4 Canada\\
{\small Email: Claude.Tardif@rmc.ca}
}

\date{July 20, 2010}

\newcommand{\cqfd}{\hfill \rule{8pt}{9pt}}
\newtheorem{define}{Definition}

\newtheorem{theorem}[define]{Theorem}

\newtheorem{remark}[define]{Remark}
\newtheorem{corollary}[define]{Corollary}

\newtheorem{conjecture}[define]{Conjecture}

\newcommand{\bd}{\begin{define} \rm}
\newcommand{\ed}{\end{define}}
\newcommand{\al}{\mbox{\rm al}}
\newcommand{\inad}{\iota}

\newcommand{\lra}{\longrightarrow}

\newfont{\Bb}{msbm10 scaled\magstep1}

\begin{document}

\maketitle
{\small \noindent{\bf Keywords:}
colourings, orientations, adjoint functors,
categorical product, Hedetniemi's conjecture,
Poljak-R\"odl function
\newline \noindent {\em AMS 2000 Subject Classification:}
05C15, 05C20, 18A40.}
\smallskip

\begin{abstract}
For an integer $k \geq 1$, the {\em $k$-th interleaved adjoint}
of a digraph $G$ is the digraph $\inad_k(G)$  with vertex-set
$V(G)^k$, and arcs $((u_1, \ldots, u_k), (v_1, \ldots, v_k))$
such that $(u_i,v_i) \in A(G)$ for $i = 1, \ldots, k$ and
$(v_i, u_{i+1}) \in A(G)$ for $i = 1, \ldots, k-1$. For every $k$ we
derive upper and lower bounds for the chromatic number of
$\inad_k(G)$ in terms of that of $G$. In the case where $G$ is
a transitive tournament, the exact value of
the chromatic number of $\inad_k(G)$ has been determined by
[H. G. Yeh, X. Zhu, 
Resource-sharing system scheduling and circular chromatic number,
Theoret. Comput. Sci. 332 (2005), 447--460].
We use the latter result in conjunction with categorial 
properties of adjoint functors
to derive the following consequence.
For every integer $\ell$, there exists a directed path $Q_{\ell}$ of algebraic length $\ell$
which admits homomorphisms into every directed graph of chromatic number at least $4$.
We discuss a possible impact of this approach on the multifactor version of
the weak Hedetniemi conjecture.
\end{abstract}

\section{Introduction}  For an integer $k \geq 1$, the {\em $k$-th interleaved adjoint}
of a digraph $G$ is the digraph $\inad_k(G)$ whose vertex-set
is $V(\inad_k(G)) = V(G)^k$, the set of all $k$-sequences of vertices of $G$,
and whose arcs-set $A(\inad_k(G))$ is the set of couples $(u,v)$ of sequences
which ``interleave'' in the sense that
for $u = (u_1, \ldots, u_k), v = (v_1, \ldots, v_k)$, we have
$(u_i,v_i) \in A(G)$ for $i = 1, \ldots, k$ and $(v_i, u_{i+1}) \in A(G)$ 
for $i = 1, \ldots, k-1$. Thus for every $k$, $\inad_k$ is a functor which 
makes a new digraph out of a given digraph $G$. We will show that interleaved adjoints share interesting
properties with iterated ``arc-digraph'' constructions $\delta^k$ (see
Section \ref{chromatic} below). Both classes of functors are categorial
right adjoints; also from the graph-theoretic
viewpoint, there are bounds for the chromatic number of $\delta^k(G)$ and $\inad_k(G)$
in terms of the chromatic number of $G$.

Here, the chromatic number $\chi(G)$ of a digraph $G$ is the minimum number of
colours needed to colour its vertices
so that pairs of vertices joined by an arc have different colours. Alternatively it is the usual
chromatic number of its symmetrisation. We view undirected graphs as symmetric digraphs.
Hence the chromatic number of a digraph $G$ is the smallest integer $n$ such that
$G$ admits a homomorphism (that is, an arc-preserving map) into $K_n$, the complete symmetric digraph
on $n$ vertices. Extending standard graph-theoretic concepts to the category
of directed graphs sometimes allow to use more elaborate categorial tools.

In Section \ref{chromatic} we find bounds for chromatic numbers of 
interleaved adjoints of general digraphs.
In particular, the case of transitive tournaments
has been dealt with in \cite{YZ}; 
in Sections \ref{adjoints}, \ref{duality} we interpret the latter result 
in terms of category theory and finite duality. 
Our main result, Theorem \ref{plh}, proves the existence of a rich class of 
paths admitting homomorphisms to all orientations of graphs with chromatic 
number at least 4.

In Section \ref{poljakrodl}, we discuss a possible connection between interleaved adjoints
and the conjecture of Hedetniemi on the chromatic number of a categorical
product of graphs: Either Theorem~\ref{plh} can be refined to the existence of families
of ``steep'' paths with similar properties, or the interleaved adjoints of transitive tournaments
witness the boundedness of a multifactor version of the (directed) ``Poljak-R\"odl'' function
from \cite{PR,TW}.

\section{Chromatic numbers of interleaved adjoints of digraphs} \label{chromatic}
The {\em arc-graph} $\delta(G)$ of a digraph $G$ is defined by
\begin{eqnarray*}
V(\delta(G)) & = & A(G)\\
A(\delta(G)) & = & \{ ((u,v),(v,w)) : (u,v), (v,w) \in A(G) \}.
\end{eqnarray*}
It is well known (see \cite{HN,PR}) that
$\log_2(\chi(G)) \leq \chi(\delta(G)) \leq 2\log_2(\chi(G))$.
Therefore for $k \geq 1$, the iterated arc graph $\delta^k(G)$ has chromatic
number in the order of $\log_2^{(k)}(\chi(G))$, where the exponential notation
represents a composition. Since $\delta^k(G)$ has no (orientations of)
odd cycles of length at most $2k+1$, the iterated arc graph construction provides
a simple constructive proof of the existence of graphs with large odd girth and
large chromatic number. The vertices of $\delta^k(G)$ correspond to chains
$(u_0, u_1, \ldots, u_k)$ of vertices of $G$ (with $(u_{i-1},u_{i}) \in A(G),
1\leq i \leq k$), and its arcs join consecutive chains
$(u_0, u_1, \ldots, u_k), (u_1, \ldots, u_k, u_{k+1})$.

We use the iterated arc graph construction to find a lower bound for the chromatic number
of interleaved adjoints of digraphs.
\begin{theorem} \label{gencol}
For an integer $k$ and a digraph $G$, we have
$$
\chi(\delta^{2k-2}(G)) \leq \chi(\inad_k(G)) \leq \chi(G).
$$
\end{theorem}

\smallskip \noindent {\em Proof.} There exists a homomorphism $\phi$ of
$\delta^{2k-2}(G)$ to $\inad_k(G)$ defined by
$$\phi(u_0, u_1, u_2, \ldots, u_{2k-2}) = (u_0, u_2, \ldots, u_{2k-2}).$$
Indeed for every arc
$((u_0, \ldots, u_{2k-2}), (u_1, \ldots, u_{2k-1}))$
of $\delta^{2k-2}(G)$, $\phi(u_0, u_1, \ldots, u_{2k-2}) =
(u_0, u_2, \ldots, u_{2k-2})$ interleaves $(u_1, u_3, \ldots, u_{2k-1})
= \phi(u_1, u_2, \ldots, u_{2k-1})$, whence $\phi$ preserves arcs.
In particular, $\chi(\delta^{2(k-1)}(G)) \leq \chi(\inad_k(G))$.

Also, there is an obvious homomorphism $\psi$ of $\inad_k(G)$ to
$G$ defined by $\psi(u_1, \ldots, u_k) = u_1$. Therefore
$\chi(\inad_k(G)) \leq \chi(G)$.
\cqfd

Even though there is a large gap between the lower and upper bounds
in Theorem~\ref{gencol}, both bounds can be tight.
The lower bound $\chi(\delta^{2k-2}(G))$ is tight when $G$
is itself of the form $\delta(H)$. The non-isolated vertices
$((u_0,u_1),(u_2,u_3), \ldots, (u_{2k-2},u_{2k-1}))$ in $V(\inad_k(G))$  then correspond
to chains  $(u_0,u_1,u_2, \ldots, u_{2k-2},u_{2k-1})$ in $V(\delta^{2k-1}(H))$,
and the arcs of $\inad_k(G)$ are interleaved sequences of the form
$$
(((u_0,u_1),(u_2,u_3), \ldots, (u_{2k-2},u_{2k-1})),
((u_1,u_2),(u_3,u_4), \ldots, (u_{2k-1},u_{2k})))
$$
which correspond to consecutive chains in $H$, that is, arcs of $\delta^{2k-1}(H)$.
Hence the arcs of $\inad_k(G)$ span a copy of
$\delta^{2k-1}(H) = \delta^{2k-2}(\delta(H)) = \delta^{2(k-1)}(G)$.
The upper bound in $\chi(G)$ is tight
in particular when $G$ is undirected. Indeed the map 
$\psi: G \lra \inad_k(G)$ where $\psi(u) = (u, u, \ldots, u)$ 
is then a homomorphism, whence the inequality 
$\chi(G) \leq \chi(\inad_k(G))$ holds. The interested reader
may find antisymmetric examples as well. 

We now turn to the case of transitive tournaments, an important class of graphs 
where $\chi(\inad_k(G))$ is linear in $\chi(G)$.
Let $T_n$ denote the transitive $n$-tournament, that is,
the vertex-set of $T_n$ is $\{ 1, \ldots, n\}$ and its arc-set is
$\{ (i,j) : 1\leq i < j \leq n \}$.

\begin{remark} \label{chi3k} For every integer $k$,
\begin{itemize}
\item[(i)] $\inad_k(T_{3k})$ contains a copy of $T_3$ induced by
$$
\{ (i, i+3, \ldots, i+3(k-1)) :  i \in \{ 1, 2, 3 \} \},
$$
\item[(ii)] $\inad_k(T_{3k})$ admits a proper 3-colouring
$c : V(\inad_k(T_{3k})) \lra \{0, 1, 2\}$ defined by
$$f(u_1, \ldots, u_k) = \left \lfloor \sum_{i=1}^k {u_i}/{k} \right \rfloor
\pmod 3.$$
\end{itemize}
Therefore $\chi(\inad_k(T_{3k})) = 3$.
\end{remark}

In \cite{YZ}, Yeh and Zhu give a general version of this result.
Let $B(n,k)$ be the symmetrisation
of $\inad_k(T_{n})$, and $K_{n/k}$ be the ``$n/k$ circular complete
graph'' in the sense of \cite{zhu2}.

\begin{theorem}[\cite{YZ}, Lemma 9] 
For integers $n \geq 2k$, 
there exist homomorphisms both ways between
$B(n,k)$ and $K_{n/k}$.
\end{theorem}

\begin{corollary} \label{chick}
For integers $n \geq 2k$,
$\chi(\inad_k(T_{n})) = \lceil n/k \rceil$. 
\end{corollary}

In fact, Yeh and Zhu proved that the circular chromatic number
of $B(n,k)$ is $n/k$, linking the ``interleaved multicolourings
of resource sharing systems'' devised by Barbosa and Gafni \cite{BG} 
to the circular chromatic number. In the next section, we examine 
the interleaved adjoints of transitive tournaments
from the point of view of category theory and finite duality.

\section{Right and Left adjoints} \label{adjoints}
The interleaved adjoints share with the arc graph construction 
the property of being categorial right adjoints. 
Following \cite{pultr,FT}, they each have
a corresponding left adjoint which acts as a kind of inverse in the
sense detailed in Theorem \ref{rladj}.
For an integer $k \geq 1$, we define the 
{\em $k$-th inverse interleaved adjoint}
of a digraph $G$ as the digraph $\iota_k^{*}(G)$ constructed
as follows.
For every vertex $u$ of $G$, we put $k$ vertices $u_1, u_2, \ldots, u_k$
in $\iota_k^{*}(G)$, and for every arc $(u,v)$ of $G$,
we put the arcs $(u_i,v_i), i = 1, \ldots, k$
and $(v_i,u_{i+1}), i = 1, \ldots, k-1$ in $\iota_k^{*}(G)$.

\begin{theorem}[\cite{pultr,FT}] \label{rladj}
For two digraphs $G$ and $H$, there exists a homomorphism
of $G$ to $\iota_k(H)$ if and only if there exists a homomorphism
of $\iota_k^{*}(G)$ to $H$.
\end{theorem}

\smallskip \noindent {\em Proof.}
If $\phi: G\lra \iota_k(H)$ is a homomorphism, we can define
a homomorphism $\psi: \iota_k^{*}(G) \lra H$ by
$\psi(u_i) = x_i$, where $\phi(u) = (x_1, \ldots, x_k)$.
Conversely, if $\psi: \iota_k^{*}(G) \lra H$ is a homomorphism,
we can define a homomorphism $\phi: G\lra \iota_k(H)$ by
$\phi(u) = (\psi(u_1), \ldots, \psi(u_k))$. 
\cqfd

\smallskip
In light of this property, we interpret the interleaved adjoints of
transitive tournaments
in term of finite path duality. For an integer $n \geq 1$, let
$P_n$ be the path with vertex-set $\{0, 1, \ldots, n\}$
and arc-set $A(P_n) = \{(0,1), (1,2), \ldots, (n-1,n)\}$.
We use the following classic result.
\begin{theorem}[\cite{gallai,hasse,roy,vitaver}] \label{classic}
A digraph $G$ admits a homomorphism to $T_n$ if and only if
there is no homomorphism of $P_n$ to $G$.
\end{theorem}
For an integer $k$, let ${\mathcal{P}}_{n,k}$ denote the
family of paths obtained from $P_n$ by reversing at most $k$ arcs.
\begin{theorem} \label{finobs}
Let $G$ be a digraph. Then the following are equivalent.
\begin{itemize}
\item[(i)] There is no homomorphism of $G$ to $\iota_k(T_n)$,
\item[(ii)] Some path in ${\mathcal{P}}_{n,k-1}$ admits a homomorphism to $G$.
\end{itemize}
\end{theorem}

\smallskip \noindent {\em Proof.} 
If there is no homomorphism of $G$ to $\iota_k(T_n)$,
then there is no homomorphism of $\iota_k^{*}(G)$ to $T_n$,
whence there exists a homomorphism $\phi: P_n \lra \iota_k^{*}(G)$.
Let $P_{\phi}$ be the path obtained from $P_n$ by reversing
the arc $(i,i+1)$ if $\phi(i) = v_j$ and $\phi(i+1) = u_{j+1}$ for some
$(u,v) \in A(G)$, and leaving it as is otherwise
(that is, if $\phi(i) = u_j$ and $\phi(i+1) = v_j$ for some
$(u,v) \in A(G)$).
Then $P_{\phi} \in {\mathcal{P}}_{n,k-1}$, and $\phi$ composed
with the natural projection of $\iota_k^{*}(G)$ to $G$ is 
a homomorphism of $P_{\phi}$ to $G$.

Conversely, suppose that some path $P$ in ${\mathcal{P}}_{n,k-1}$
admits a homomorphism $\phi: P \lra G$. Let $f: V(P) \lra \mbox{\Bb Z}$
be the function defined recursively by $f(0) = 1$ and
$$
f(i+1) = \left \{
\begin{array}{l}
\mbox{$f(i)$  if $(i,i+1) \in A(P)$}, \\
\mbox{$f(i) + 1$  if $(i+1,i) \in A(P)$}.
\end{array}
\right.
$$
Since $P$ is in ${\mathcal{P}}_{n,k-1}$, we have $1 \leq f(i) \leq k$ for all 
$i$.
If $(i,i+1) \in A(P)$, then $(\phi(i),\phi(i+1)) \in A(G)$ and 
$f(i+1) = f(i)$,
hence $(\phi(i)_{f(i)},\phi(i+1)_{f(i+1)}) \in A(\iota_k^{*}(G))$.
If $(i+1,i) \in A(P)$, then $(\phi(i+1),\phi(i)) \in A(G)$ and 
$f(i+1) = f(i)+1$,
hence $(\phi(i)_{f(i)},\phi(i+1)_{f(i+1)}) \in A(\iota_k^{*}(G))$.
Therefore the map $\psi: P_n \lra \iota_k^{*}(G)$
defined by $\psi(i) = \phi(i)_{f(i)}$ is a homomorphism.
This implies that there is no homomorphism of $\iota_k^{*}(G)$
to $T_n$, and no homomorphism of $G$ to $\iota_k(T_n)$.
\cqfd

Both \cite{BG} and \cite{YZ} provide polynomial algorithms to decide whether
an input digraph $G$ admits a homomorphism to $\iota_k(T_n)$. From the point
of view of descriptive complexity, Theorem~\ref{finobs} shows that the problem
has ``width 1'' in the sense of \cite{FV}, and is solvable polynomially through the
``arc consistency check'' algorithm. More precisely, $\iota_k(T_n)$ has ``finite duality''.
We will discuss the structural consequences of this fact in the next section.
We end this section by noting the following consequence of Theorem~\ref{finobs} and
Corollary \ref{chick}.
\begin{corollary} \label{minty}
For every integers $c, k$ and for every digraph $G$
such that $\chi(G) > c$, there exists a path 
$P \in {\mathcal{P}}_{ck,k-1}$ which admits a homomorphism to $G$. \cqfd
\end{corollary}
Actually, this is also an easy consequence of ``Minty's painting lemma''~\cite{minty}.
We include it here for reference; we will discuss hypothetical
strenghtenings of it in Section~\ref{poljakrodl}.

\section{Finite duality} \label{duality}
A digraph $H$ has {\em finite duality} if there exists a finite
family $\mathcal{F}$ of digraphs such that a graph $G$ admits a
homomorphism to $H$ if and only if there is no homomorphism
of a member of $\mathcal{F}$ to $G$. The family $\mathcal{F}$
is then called a {\em complete set of obstructions} for $H$.
For instance, by Theorem~\ref{classic},
$T_n$ has finite duality and admits $\{P_n\}$ as
a complete set of obstructions. By Theorem~\ref{finobs},
$\inad_k(T_n)$ also has finite duality and admits ${\mathcal{P}}_{n,k-1}$
as a complete set of obstructions.

The finite dualities were characterised in \cite{NT1}, in terms
of homomorphic equivalence with categorical products of structures
with singleton duality.
Two digraphs $G$, $H$ are called {\em homomorphically equivalent}
if there exist homomorphisms of $G$ to $H$ and of $H$ to $G$.
The {\em categorical product} of a family
$\{G_1, \ldots, G_n\}$ of digraph is the digraph $\Pi_{i=1}^n G_i$ defined by
\begin{eqnarray*}
V(\Pi_{i=1}^n G_i) & = & \Pi_{i=1}^n V(G_i),\\
A(\Pi_{i=1}^n G_i) & = & \{ ((u_1,\ldots,u_n),(v_1,\ldots,v_n)) : (u_i,v_i) \in
A(G_i)
\mbox{ for } 1 \leq i \leq n \}.
\end{eqnarray*}
We use mostly categorical products of sets of digraphs.
For ${\mathcal F} = \{G_1, \ldots, G_n\}$, we write $\Pi {\mathcal F}$
for $\Pi_{i=1}^n G_i$. This allows to simplify the notation without
loss of generality, since the categorical product is commutative and
associative (up to isomorphism).

\begin{theorem}[\cite{NT1}] \label{findu}
For every directed tree $T$, there exists a directed graph $D(T)$
(called the {\em dual} of $T$)
which admits $\{ T \}$ as a complete set of obstructions.
A digraph $H$ has finite duality
if and only if there exists a finite family
${\mathcal{F}}$ of trees such that $H$ is homomorphically equivalent
to $\Pi \{D(T) : T \in {\mathcal{F}}\}$. ${\mathcal{F}}$ is then a complete
set of obstructions for $H$.
\end{theorem}
\begin{corollary} \label{inadprod}
$\inad_k(T_n)$ is homomorphically equivalent to 
$\Pi \{D(P) : P \in {\mathcal{P}}_{n,k-1}\}$.
\end{corollary}

The dual $D(T)$ of a tree $T$ is not unique, but all duals
of a given tree are homomorphically equivalent.
According to \cite{NT2}, we get a possible construction for $D(T)$
by taking for $V(D(T))$ the set of all functions $f: V(T) \lra A(T)$
such that $f(u)$ is incident to $u$ for every $u \in V(T)$,
and putting an arc from $f$ to $g$ if for every $(u,v) \in A(T)$,
$f(u) \neq g(v)$. This is the simplest general construction known,
and yet is far from transparent. This combined with the fact that
the family ${\mathcal{P}}_{n,k-1}$ is large makes it difficult to
describe explicit homomorphisms between $\inad_k(T_n)$ and
$\Pi \{D(P) : P \in {\mathcal{P}}_{n,k-1}\}$.
However, combined with the concept of algebraic length, we use
this structural insight to exhibit some common
features of digraphs with a large chromatic number.

The {\em algebraic length}  of a path $P$ is the value
$$
\al(P) = \min \{ n : \mbox{there exists a homomorphism of
$P$ to $P_n$} \}.
$$
If we picture a path drawn from left to right, 
then its algebraic length is the
(absolute value of) the difference between the number of its 
forward arcs and
the number of its backward arcs. In particular, for 
$n \geq k \geq 1$, the paths in
${\mathcal{P}}_{n,k-1}$ all have algebraic length at least
$n - 2k + 2$.

We use the following results.
\begin{theorem}[\cite{HZ}] \label{hompath}
A digraph $G$ admits a homomorphism to $P_n$ if and only if no path
of algebraic length $n+1$ admits a homomorphism to $G$.
\end{theorem}
\begin{theorem}[\cite{NP}] \label{mulpath}
A categorical product $\Pi_{i=1}^n G_i$ 
of digraphs admits a homomorphism to $P_n$
if and only if at least one of the factors does.
\end{theorem}

Our main result is the following.
\begin{theorem} \label{plh}
For every length $\ell$,
there exists a path $Q_{\ell}$ such that $\al(Q_{\ell}) = \ell$ and for every
digraph $G$ with chromatic number at least $4$, there exists a homomorphism
of $Q_{\ell}$ to $G$.
\end{theorem}

\smallskip \noindent {\em Proof.} 
For $\ell \leq 3$, we can take $Q_{\ell} = P_{\ell}$
by Theorem~\ref{classic}. For $\ell \geq 4$, we put $k = \ell-2$.
The paths in ${\mathcal{P}}_{3k,k-1}$ all have algebraic length
at least $\ell$, hence none of them admits a homomorphism to
$P_{\ell-1}$.  Thus by Theorem~\ref{mulpath}, their categorical product
$\Pi {{\mathcal{P}}_{3k,k-1}}$ does not admit a homomorphism to
$P_{\ell-1}$. Therefore by Theorem~\ref{hompath}, there exists a path
$Q_{\ell}$ of algebraic length $\ell$ which admits a homomorphism
to $\Pi {{\mathcal{P}}_{3k,k-1}}$. We will show that
$Q_{\ell}$ has the required property.

Since $Q_{\ell}$ admits a homomorphism to 
$\Pi {{\mathcal{P}}_{3k,k-1}}$,
it admits a homomorphism to every path 
$P$ in ${\mathcal{P}}_{3k,k-1}$.
By Theorem~\ref{findu}, this implies that there 
is no homomorphism of
$P$ to $D(Q_{\ell})$. Since this holds for every 
$P$ in ${\mathcal{P}}_{3k,k-1}$,
by Theorem~\ref{finobs} there exists a homomorphism of
$D(Q_{\ell})$ to $\inad_k(T_{3k})$. Therefore by Lemma~\ref{chick},
we have $\chi(D(Q_{\ell})) \leq \chi(\inad_k(T_{3k})) = 3$.

Now let $G$ be a digraph such that
there is no homomorphism of $Q_{\ell}$ to $G$. 
Then by Theorem~\ref{findu}
there exists a homomorphism of $G$ to $D(Q_{\ell})$ whence
$\chi(G) \leq \chi(D(Q_{\ell})) \leq 3$.
Therefore $Q_{\ell}$ admits homomorphisms 
to all digraphs with chromatic
number at least 4. 
\cqfd

The study of the relation between the algebraic length of paths
and the chromatic number of their duals was 
initiated in~\cite{NT2}, where it was shown that the bound
$\chi(D(P)) < \al(P)$ can hold for paths with arbitrarily large
algebraic lengths. The paper~\cite{NT3} provides examples 
of paths $P$ such that $\chi(D(P)) \in O(\log(\al(P)))$,
and raises the question of the existence of paths with arbitrarily
large algebraic lengths whose duals have bounded chromatic number.
Theorem~\ref{plh} settles this question in the affirmative.
The focus now shifts to how steep can such paths be,
that is, how many arcs must there be in a path with large
algebraic length whose dual has a small chromatic number.
As we shall see, this question is connected to a long-standing
conjecture in graph theory.

\section{The multifactor Poljak-R\"odl function} \label{poljakrodl} 

The Poljak-R\"odl function $f: \mbox{\Bb N} \lra \mbox{\Bb N}$
for directed graphs is defined by
\begin{eqnarray*}
f(c) & = & \min \{ \chi(G_1 \times G_2) :
\mbox{$G_1$ and $G_2$ are $c$-chromatic digraphs} \}.
\end{eqnarray*}
Its undirected version $g: \mbox{\Bb N} \lra \mbox{\Bb N}$
is defined by
\begin{eqnarray*}
g(c) & = & \min \{ \chi(G_1 \times G_2) :
\mbox{$G_1$ and $G_2$ are $c$-chromatic undirected graphs} \}.
\end{eqnarray*}
 These functions are related to the long-standing
conjecture of Hedetniemi, first formulated in 1966:
\begin{conjecture}[\cite{hedetniemi}] \label{hedcon}
If $G$ and $H$ are undirected graphs, then
$$\chi(G\times H) = \min \{ \chi(G), \chi(H)\}.$$
\end{conjecture}
Hedetniemi's conjecture states that $g(c) = c$
for all $c$, but for the moment it is not even known
whether $g$ is bounded or unbounded. 
In \cite{TW}, it is shown that $f$ is unbounded if and 
only if $g$ is unbounded, and in \cite{poljak,zhu},
it is shown that $f$ is bounded if and only if 
$f(c) \leq 3$ for all $c$. 

We will now consider an application of this topic to a variation
of Corollary~\ref{minty}.
\begin{conjecture} \label{minty2}
There exists a number $C$ such that 
for every $k$ and every partition
of ${\mathcal{P}_{3k,k-1}}$ in two sets ${\mathcal{Q}_{1}}$,
${\mathcal{Q}_{2}}$, we can select a set ${\mathcal{Q}_{i}}$
with the following property:
\begin{quote}
For every digraph $G$ such that $\chi(G) \geq C$, there exists
a path $P \in {\mathcal{Q}_{i}}$ which admits a homomorphism
to $G$.
\end{quote}
\end{conjecture}

\begin{theorem}
If $f$ is unbounded, then Conjecture~\ref{minty2} is true.
\end{theorem}

\smallskip \noindent {\em Proof.} Suppose that
$f$ is unbounded, and let $C$ be the first value
such that $f(C) > 3$. When
${\mathcal{P}_{3k,k-1}}$ is partitioned in two sets
${\mathcal{Q}_{1}}$, ${\mathcal{Q}_{2}}$,
we have
$$
\Pi \{ D(P) : P \in {\mathcal{P}_{3k,k-1}} \}
=
\left ( \Pi \{ D(P) : P \in  {\mathcal{Q}_{1}} \} \right )
\times
\left ( \Pi \{ D(P) : P \in  {\mathcal{Q}_{2}} \} \right ).
$$
We have $\chi(\Pi \{ D(P) : P \in {\mathcal{P}_{3k,k-1}} \})
= \chi(\inad_k(T_{3k})) = 3$
by Remark~\ref{chi3k} and Corollary~\ref{inadprod}.
Since $f(C) > 3$, one of the factors
$\Pi \{ D(P) : P \in  {\mathcal{Q}_{i}} \}$ has chromatic 
number less than $C$. Therefore for every digraph $G$ such that 
$\chi(G) \geq C$, there exists no homomorphism of $G$
to $\Pi \{ D(P) : P \in  {\mathcal{Q}_{i}} \}$, whence
there exists
a path $P \in {\mathcal{Q}_{i}}$ which admits a homomorphism
to $G$. \cqfd

We now turn to multifactor versions of the Poljak-R\"odl functions.
Define  $f_{\rm m}, g_{\rm m}: \mbox{\Bb N} \lra \mbox{\Bb N}$ by
\begin{eqnarray*}
f_{\rm m}(c) & = & \min \{ \chi(\Pi {\mathcal{D}}) :
\mbox{${\mathcal{D}}$ is a finite family of $c$-chromatic digraphs} \} \\
g_{\rm m}(c) & = & \min \{ \chi(\Pi {\mathcal{G}}) :
\mbox{${\mathcal{G}}$ is a finite family of $c$-chromatic graphs} \}.
\end{eqnarray*}
Hedetniemi's conjecture is obviously equivalent 
to its multifactor version, that is, the identity
$\chi(\Pi {\mathcal{G}}) = \min \{ \chi(G) : G \in {\mathcal{G}} \}$
for every finite family ${\mathcal{G}}$ of undirected graphs.
(Miller \cite{miller} notes that the identity can fail if
${\mathcal{G}}$ is infinite.) However it is not known whether
the unboundedness of $g_{\rm m}$ would imply that of $f_{\rm m}$.
These hypotheses involving
various versions of the Poljak-R\"odl functions are summarized
in the following table.
$$
\begin{array}{ccc}
\mbox{$g(c) \equiv c$} & \Leftrightarrow & \mbox{$g_{\rm m}(c) \equiv c$} \\
\Downarrow & & \Downarrow \\
\mbox{$g(c)$ is unbounded} & \Leftarrow & \mbox{$g_{\rm m}(c)$ is unbounded} \\
\Updownarrow & & \Uparrow \\
\mbox{$f(c)$ is unbounded} & \Leftarrow & \mbox{$f_{\rm m}(c)$ is unbounded}
\end{array}
$$
At the top of the diagram is Hedetniemi's conjecture, which has withstood
scrutiny for more than forty years. At the bottom right is the hypothesis of
the unboundedness of $f_{\rm m}$, which is very much in doubt. Consider the function
$h: \mbox{\Bb N} \lra \mbox{\Bb N}$ defined by
\begin{eqnarray*}
h(k) & = & \min \{ \chi(D(P)) : P \in {\mathcal{P}_{3k,k-1}}\}.
\end{eqnarray*}
If $h$ is unbounded, then the products
$\Pi \{ D(P) : P \in {\mathcal{P}_{3k,k-1}} \}$ witness
the boundedness of $f_{\rm m}$. On the other hand,
if $h(k) < C$ for all $k$, then every set
${\mathcal{P}_{3k,k-1}}$ contains a path $P$ such that $\chi(D(P)) < C$.
This would imply Conjecture~\ref{minty2}, as well as a variation of
of Theorem~\ref{plh} for chromatic nomber $C$, using steep paths rather than 
the family $\{\mathcal{Q}_{\ell} \}_{\ell \in \mbox{\Bb N}}$. 
All in all, it seems
to be too strong a certificate for $(C-1)$-colourability.

It would be interesting to know whether a family 
$\{Q_{\ell}\}_{\ell \in \mbox{\Bb N}}$ of paths 
which respects the conclusion of Theorem~\ref{plh}
needs an exponential number of arcs;
and also which of the missing implications can be added
to the table above.

\smallskip \noindent {\bf Acknowledgements.} The results \cite{YZ} 
of Yeh and Zhu were communicated to the authors by Reza Naserasr
during the Fields Workshop on Cocliques and Colourings
at the University of Waterloo (June 1-5, 2009).

\end{document}